\documentclass{amsart}
\usepackage{amsmath,amssymb}

\def\ind{{\mathcal I}}
\def\res{{\mathcal R}}
\providecommand{\sgn}{\mathop{\rm sgn}\nolimits}
\providecommand{\car}{\mathop{\rm char}\nolimits}
\providecommand{\End}{\mathop{\rm End}\nolimits}

\newtheorem{Lemma}{Lemma}[section]
\newtheorem{Theorem}[Lemma]{Theorem}
\newtheorem{Corollary}[Lemma]{Corollary}
\newtheorem{Proposition}[Lemma]{Proposition}

\begin{document}

\title{Branching Rules for Specht Modules}


\author{Harald Ellers}
\address{Department of Mathematics\\
   Northern Illinois University\\
   DeKalb, IL  60115\\
   USA}
\email{ellers@math.niu.edu}
\author{John Murray}
\address{Department of Mathematics\\
   National University of Ireland - Maynooth\\
   Co.~Kildare\\
   Ireland}
\email{jmurray@maths.may.ie}

\subjclass{20C20, 20C30}
\date{May 20, 2004}

\begin{abstract}
Let $\Sigma_n$ be the symmetric group of degree $n$, and let $F$ be a field of characteristic
distinct from $2$. Let $S_F^{\lambda}$ be the Specht module over $F\Sigma_n$ corresponding to the partition $\lambda$ of $n$. We find the indecomposable components of the restricted module
$S_F^{\lambda}\downarrow_{\Sigma_{n-1}}$ and the induced module 
$S_F^{\lambda}\uparrow^{\Sigma_{n+1}}$. Namely, if $b$ and $B$ are block idempotents 
of $F\Sigma_{n-1}$ and $F\Sigma_{n+1}$ respectively, then the modules 
$S_F^{\lambda}\downarrow_{\Sigma_{n-1}}b$ and 
$S_F^{\lambda}\uparrow^{\Sigma_{n+1}}B$ are $0$ or indecomposable. We give examples
to show that the assumption $\car F \neq 2$ cannot be dropped.
\end{abstract}

\maketitle

\section{Introduction}

Let $n$ be a positive integer and let $\Sigma_n$ be the 
symmetric group of degree $n$. 
For any field $F$ and any partition $\lambda$ of $n$, 
the Specht module $S_{F}^{\lambda}$ is defined to be the submodule
of the permutation module $F_{\Sigma_{\lambda}}\uparrow ^{\Sigma_n}$ 
spanned by certain elements called polytabloids, where $\Sigma_{\lambda}$ is 
the Young subgroup associated to $\lambda$
and $F_{\Sigma_{\lambda}}$ is the principal 
$ F\Sigma_{\lambda}$-module. (See \cite{James}
for definitions.)
Specht modules play a central role in the representation theory
of the symmetric group, because in characteristic $0$ 
the Specht modules are the simple $F\Sigma_n$-modules, while  
in characteristic $p$ the heads of the Specht 
modules   $S_{F}^{\lambda}$ such that  $\lambda$ is  $p$-regular 
are the simple $F\Sigma_n$-modules.
When the field $F$ has charactersitic $0$, the structure of 
the restriction of $S_{F}^{\lambda}$ to $\Sigma_{n-1}$ is given 
by the Classical Branching Rule: the module 
$S_{F}^{\lambda}\downarrow_{\Sigma_{n-1}}$ is 
a direct sum $\bigoplus_{\mu}S_{F}^{\mu}$, where $\mu$ runs through all 
partitions of $n-1$ obtained from $\lambda$ by removing a node 
from its Young diagram. In 1971, Peel \cite{Peel} gave the first 
characteristic $p$ version of the branching rule. He showed that 
there is a series of submodules such that 
the successive quotients are the Specht modules $S_{F}^{\mu}$, 
where $\mu$ runs through the same set. Nevertheless, the 
structure of the restriction $S_{F}^{\lambda}\downarrow_{\Sigma_{n-1}}$
is not well understood. For example, the problem of finding a 
composition series is open and very difficult, 
and the socle is not known. See Kleshchev \cite{Kl}
for an introduction to recent work on 
$S_{F}^{\lambda}\downarrow_{\Sigma_{n-1}}$.

In this paper, we find the indecomposable components of 
$S_{F}^{\lambda}\downarrow_{\Sigma_{n-1}}$, 
when the characteristic of $F$
is not $2$.  These are given by Theorem \ref{Main}: 
if $b$ is a block idempotent of $F\Sigma_{n-1}$, then 
$S_{F}^{\lambda}\downarrow_{\Sigma_{n-1}}b$ 
is $0$ or indecomposable.  Thus there is a bijection between 
the set of   indecomposable components 
of  $S_{F}^{\lambda}\downarrow_{\Sigma_{n-1}}$  and 
the set of $p$-cores that can be obtained from $\lambda$ by removing
first one node and then a sequence of rim $p$-hooks. 
We also prove the analogous theorem for the induced module
$S_{F}^{\lambda}\uparrow^{\Sigma_{n+1}}$.  The two proofs are 
almost identical. 
We give examples to 
show that the assumption $\car F \neq 2$ cannot be dropped.

The combinatorial part of the proof is in section 2. Here we find the 
minimal polynomials for the actions of $E_{n-1}$ on 
 $S_{F}^{\lambda}\downarrow_{\Sigma_{n-1}}$ and 
 $E_{n+1}$ on  $S_{F}^{\lambda}\uparrow^{\Sigma_{n+1}}$,
 where  $E_k$ is the sum of all the transpositions in $\Sigma_k$. 
 These polynomials have degrees $m$ and $m+1$ respectively, 
 where $m$ is the number of distinct
 parts of $\lambda$.  The results of section 2 are valid for all fields, 
 not just those of odd characteristic. 
 
 In section 3, we investigate the algebras ${\mathcal E}= 
 \End_{F\Sigma_{n-1}}(S_{F}^{\lambda}\downarrow_{\Sigma_{n-1}})$
 and ${\mathcal F} = \End_{F\Sigma_{n+1}}(S_{F}^{\lambda}\uparrow^{\Sigma_{n+1}})$. 
  Under the assumption that  $\car F \neq 2$,
 we use the results from section 2 to show that the natural maps
 $ Z(F\Sigma_{n-1}) \rightarrow {\mathcal E}/J({\mathcal E})$ and 
 $ Z(F\Sigma_{n+1}) \rightarrow {\mathcal F}/J({\mathcal F})$ are surjective,
where $J({\mathcal E})$  and  $J({\mathcal F})$ are the Jacobson radicals
of  $\mathcal E$ and $\mathcal F$. 
 The main theorem follows easily.

\section{The minimal polynomials of the sum of all 
          transpositions acting on
          the restriction and induction of a Specht module}

Throughout this paper $n$ is a fixed positive integer and
$\lambda$ is a fixed partition of $n$. We orient the Young
diagram $[\lambda]$ left to right and top to bottom.
This means that 
the  first row is the one at the top and the first
column is the one at the left. The $(i,j)$ node is in the $i$th
row and the $j$th column.
We will use ${\widehat n}$ to denote the set $\{1,\ldots,n\}$
and let $\Sigma_n$ denote the group of permutations of ${\widehat n}$.
Permutations and homomorphisms will generally act on the right.
The \emph{Murphy element} $L_n$ is the sum of all transpositions
in $\Sigma_n$ that are not in $\Sigma_{n-1}$ (with $L_1:=0$). We use
$E_n$ to denote the sum of all transpositions in $\Sigma_n$. So $E_n$
is the $1$-st elementary symmetric function in the Murphy elements.

Let $F$ be any field and let $S^{\lambda}$ denote the Specht
module, defined over $F$, corresponding to $\lambda$.
We use the notation
\[
\begin{aligned}{}
\res\quad&\mbox{for the restriction of $S^\lambda$ to $\Sigma_{n-1}$ 
and}\\
\ind\quad&\mbox{for the induction of $S^\lambda$ to $\Sigma_{n+1}$}.
\end{aligned}
\]
The purpose of this section is to compute
the minimal polynomial of $E_{n-1}$ acting on $\res$ and
the minimal polynomial of $E_{n+1}$ acting on $\ind$.

We consider a \emph{$\lambda$-tableau} to be a bijective map
$t:[\lambda]\rightarrow{\widehat n}$. The value of $t$ at a node
$(r,c)$ is denoted by $t_{rc}$. The group $\Sigma_n$ acts on
the set of all 
$\lambda$-tableaux by functional composition;
$(t\pi)_{rc}=t_{rc}\pi$, for each $\pi\in\Sigma_n$.

Suppose that $\lambda$ has $l$ nonzero parts
$[\lambda_1\geq\lambda_2\geq\ldots\geq\lambda_l]$.
We regard a \emph{$\lambda$-tabloid} as an ordered partition
${\mathcal P}=({\mathcal P}_1,\ldots,{\mathcal P}_l)$ of ${\widehat n}$ such that
the cardinality of ${\mathcal P}_u$ is $\lambda_u$, for $u=1,\ldots,l$.
Each $\lambda$-tableau $t$ determines the $\lambda$-tabloid $\{t\}$
whose $u$-th part is the set of entries in the $u$-th row of $t$.
If $s$ is a $\lambda$-tableau, then $\{t\}=\{s\}$ if and only if
$s=t\pi$, for some $\pi$ in the row stabilizer $R_t$ of $t$.
We denote the column stabilizer of $t$ by $C_t$.
We denote by $M^{\lambda}$ the $F\Sigma_n$-module
consisting of  all formal $F$-linear combinations of $\lambda$-tabloids.

Adapting the notation of James \cite{James}, let 
$(r_1,c_1),\ldots,(r_m,c_m)$
be the removable nodes of $[\lambda]$, ordered so that
$r_1<\ldots<r_m$ and $c_1>\ldots>c_m$. Set $r_0=0=c_{m+1}$.
The addable nodes of $[\lambda]$ are the $(m+1)$ nodes
$(r_u+1,c_{u+1}+1)$, for $u=0,\ldots,m$.
We use $\lambda{\downarrow_u}$ to denote the partition of $n-1$
obtained by decrementing the $r_u$-th part of $\lambda$ by $1$,
for $u\in\widehat{m}$. In addition, we use $\lambda{\uparrow}^u$ to
denote the partition of $n+1$ obtained by incrementing the
$(r_u+1)$-th part of $\lambda$ by $1$, for $u\in\widehat{m+1}$.

We need special notation for certain subsets of 
a $\lambda$-tableau $t$.
For the rest of the paper, suppose that
$\lambda$ has parts of $m$ different nonzero lengths.
For any $u \in {\widehat m}$,
let $H_u(t)$ be the set of entries in the union of the top $r_u$
rows of $t$, and let $V_u(t)$ be the set of entries in the union
of columns of $t$ numbered from $c_{u+1}+1$ to $c_{u}$ (inclusive).
Clearly $H_1(t)\subset\ldots\subset H_m(t)$, while
$V_m(t),\ldots,V_1(t)$ forms a partition of $t$. Also
$V_u(t)\subseteq H_v(t)$ if and only if $u\leq v$.
As $H_u(t)$ depends only on the rows of $t$,
we may define $H_u(\{t\}):=H(t)$.

By Theorem 9.3 in \cite{James}, $\res$ has a Specht series
\[
0\subset\res_1\subset\res_2\subset\ldots\subset\res_{m}=\res,
\]
with $\res_u/\res_{u-1}\cong S^{\lambda{\downarrow_u}}$,
for $u\in{\widehat m}$. Also, by 17.14 in \cite{James},
$\ind$ has a Specht series
\[
\ind=\ind_1\supset\ind_2\supset\ldots\supset\ind_{m+1}\supset\ind_{m+2}=0,
\]
with $\ind_u/\ind_{u+1}\cong S^{\lambda{\uparrow^u}}$,
for $u\in{\widehat{m+1}}$. Each factor $\ind/\ind_{u+1}$
is isomorphic to a submodule of the permutation module
$M^{\lambda{\uparrow^u}}$.

\begin{Lemma}\label{L:central}
Suppose that the $F\Sigma_n$-module $M$ has a Specht
series $0=M_0\subset M_1\subset\ldots\subset M_m=M$.
Let $z\in Z(F\Sigma_n)$ and let $u\in{\widehat m}$. Then
there is a scalar $z_u$ in $F$ such that the map
$M_u/M_{u-1}\rightarrow M_u/M_{u-1}$ given by multiplication
by $z$ is equal to $z_u$ times the identity map.
\end{Lemma}

\begin{proof}
If $\car F=0$, then $M_u/M_{u-1}$ is an irreducible $F\Sigma_n$-module
(a Specht module), and the conclusion is obvious.
If $\car F =p$ is positive, then $M_u/M_{u-1}$ is the $p$-modular reduction
of an irreducible module defined over a suitable discrete valuation ring
of characteristic $0$. The conclusion follows in this case from the
characteristic zero case. 
\end{proof}

This lemma allows us to give the following upper bound on the
degrees of the minimal polynomials of $E_{n-1}$ and $E_{n+1}$.

\begin{Corollary}\label{C:min-poly-upper-bound}
The minimal polynomial of $E_{n-1}$ acting on $\res$
has degree at most $m$, while the minimal polynomial
of $E_{n+1}$ acting on $\ind$ has degree at most $m+1$.
\end{Corollary}

\begin{proof}
Let $u\in{\widehat m}$. Lemma \ref{L:central}
shows that $\res_u(E_{n-1}-z_u)\subseteq \res_{u-1}$,
for some scalar $z_u$. It follows from a simple inductive argument
that $\res\prod_{u=1}^m(E_{n-1}-z_u)=0$. A similar argument deals
with the action of $E_{n+1}$ on $\ind$. 
\end{proof}

It will turn out that the polynomials given in the proof
of Corollary \ref{C:min-poly-upper-bound} are minimal.
Before we prove this, we will identify
the scalars $z_u$ in terms of Young diagrams.

The residue of a node $(r,c)$ is the scalar $(c-r)1_F$.
We set $E(\lambda)$ as the sum of the residues of all nodes in $[\lambda]$.
So $E(\lambda)$ is the $1$-st elementary symmetric function in the residues.
An easy calculation shows that
$E(\lambda)=\sum_{i=1}^l\frac{1}{2}\lambda_i(\lambda_i+1-2i)1_F$.
The next lemma is a special case of a more general result
proved by G. E. Murphy \cite{Murphy}:
$1$-st elementary symmetric function can be
replaced by any symmetric function in $n$ variables.

\begin{Lemma}\label{L:En-action}
$E_n$ acts as the scalar $E(\lambda)$ on $S^\lambda$.
\end{Lemma}

\begin{proof}
Let $t$ be a $\lambda$-tableau, let $(r,c)\in[\lambda]$
and let $i=t_{rc}$. Fix $1\leq c_{\prime}<c$. Then by a simple
Garnir relation (section 7 of \cite{James}),
$e_t\sum_{j}(i,j)=e_t$, where $j$ runs
over all entries in the $c_{\prime}$-th column of $t$.
Also $e_t(i,j)=-e_t$, for each entry $j$ above $i$ in
column $c$ of $t$. It follows that
\[
e_t\sum_{j}(i,j)=(c-r)e_t,
\]
where $j$ runs over those elements of ${\widehat n}$
that lie in $t$ in columns strictly left of $i$
or in the same column as $i$ but strictly 
above $i$.
If we sum over all $(r,c)\in[\lambda]$, each transposition $(i,j)$
occurs exactly once on the left hand side, while the coefficient of
$e_t$ on the right hand side is $E(\lambda)$.
\end{proof}

If $t$ is a $\lambda$-tableau, the polytabloid
$e_t$ is the following element of $M^\lambda$:
\[
e_t:=\sum_{\pi\in C_t}\sgn \pi \{t\pi\}.
\]
It is well known that the polytabloids span the Specht module $S^\lambda$.
James' description of $\res$, and the Garnir relations, show that $e_t$
lies in $\res_u\backslash\res_{u-1}$ if $n\in V_u(t)\backslash H_{u-1}(t)$
(although we do not use this fact).

We next describe the induced module $\ind$.
Suppose that $u\in{\widehat{m+1}}$. Let $T$ be a
$\lambda{\uparrow^u}$-tableau, and let $t$ denote
the restriction of $T$ to $[\lambda]$.
Then the $(\lambda,T)$-polytabloid $e_T^\lambda$
is the following element of $M^{\lambda{\uparrow^u}}$:
\[
e_T^\lambda:=\sum_{\pi\in C_t} \sgn \pi \{T\pi\}.
\]
In Section 17 of \cite{James}, James
has shown that when $u=m+1$, the corresponding
$(\lambda,T)$-polytabloids span
an  $F\Sigma_{n+1}$-submodule of $M^{\lambda{\uparrow^{m+1}}}$,
which is isomorphic to the induced module $\ind$. We will always
work with this copy of $\ind$.

When we are showing that the polynomials given in the
proof of \ref{C:min-poly-upper-bound} are minimal, it will
be convenient to look at the action of the Murphy
elements $L_{n}$ and $L_{n+1}$ rather than $E_{n-1}$ and
$E_{n+1}$. The following lemma provides a link between
these actions. If $t$ is a $\lambda$-tableau, 
its {\em extension to } $[\lambda{\uparrow^{m+1}}]$
is the $\lambda{\uparrow^{m+1}}$-tableau that is obtained
from $t$ by appending $n+1$ to the bottom of the
first column. 

\begin{Lemma}\label{L:poly}
Let $t$ be a $\lambda$-tableau and let $T$ be its extension to
$[\lambda{\uparrow^{m+1}}]$. Suppose that $f(x)\in F[x]$. Then
\[
\begin{aligned}{}
e_t         \,f(E_{n-1})\quad&=&e_t         \,&f(E(\lambda)-L_n);\\
e_T^\lambda \,f(E_{n+1})\quad&=&e_T^\lambda 
\,&f(E(\lambda)+L_{n+1}).
\end{aligned}
\]
\end{Lemma}

\begin{proof}
Lemma \ref{L:En-action} shows that $E_n$ acts as the scalar $E(\lambda)$
on $\res$. The first statement then follows from $E_{n-1}=E_n-L_n$.

Consider the subspace $V$ of $M^{\lambda{\uparrow^{m+1}}}$ spanned by
all $e_U^\lambda$ such that $U$ is a
$\lambda{\uparrow^{m+1}}$-tableau with $n+1$ in the unique
entry of its last row. The subspace $V$ is a
direct summand of  the restriction of
$\ind$ to $\Sigma_n$, and as an $F\Sigma_n$-submodule,
$V$ is clearly isomorphic to $S^\lambda$.
Thus $e_T^\lambda$ lies in a direct summand of the restriction of
$\ind$ to $\Sigma_n$ that is isomorphic to $S^\lambda$. So Lemma
\ref{L:En-action} shows that $e_T^\lambda E_n=E(\lambda)e_T^\lambda$.
The second statement now follows from $E_{n+1}=E_n+L_{n+1}$,
and the fact that $E_nL_{n+1}=L_{n+1}E_n$. 
\end{proof}

When we are showing that the polynomials given in the
proof of \ref{C:min-poly-upper-bound} are minimal,
we will want to show that there is a $\lambda$-tableau $t$
such that the set of vectors
$\{e_t (L_n)^i\mid 0\leq i\leq m-1\}$ is linearly independent.
This will be accomplished using the following
technical lemma concerning the action of $L_n$ on $\res$.

\begin{Lemma}\label{L:L_n^i}
Let $t$ be a $\lambda$-tableau such that $n\in V_m(t)\backslash H_{m-1}(t)$.
For each $u\in \widehat{m-1}$, choose $x_u\in V_u(t)\backslash 
H_{u-1}(t)$.
Set $s=t\,(n,x_{m-1},x_{m-2},\ldots,x_1)$.
Let $i$ be a positive integer with $i \leq m-1$.
Then the coefficient  of $\{s\}$ in the expansion
of $e_t (L_n)^i$ into tabloids is
\[
\begin{aligned}{}
&0,\quad\mbox{when $0\leq i\leq m-2$;}\\
&1,\quad\mbox{when $i=m-1$.}
\end{aligned}
\]
\end{Lemma}

\begin{proof}
Clearly $(L_n)^i=\sum(w_i,n)(w_{i-1},n)\ldots(w_1,n)$, where
$(w_1,\ldots,w_i)$ ranges over all functions
${\widehat i}\rightarrow{\widehat{n-1}}$.
Let $(y_1, \ldots, y_{i})$ be a function
${\widehat i}\rightarrow{\widehat{n-1}}$, let
$\theta=(y_i,n)(y_{i-1},n)\ldots(y_1,n)$,
  and assume that
$\{s\}$ appears with nonzero coefficient in the expansion of
$e_t \theta$. We have two goals: (a) 
to show that $i = m-1$, and (b)
to show that when $i=m-1$, the sequence $(y_1, \ldots, y_{m-1})$
is equal to the sequence $(x_1, \ldots, x_{m-1})$.
The second part of the lemma follows easily from this
second goal, as we now show.  In the sum
$\sum e_t (w_i n) \ldots (w_1 n)$, $\{ s \}$ can appear
in only one term, namely $e_t (x_{m-1}, n) \ldots (x_{1}, n)$.
Since this term is equal to
$e_{t}(n,x_{m-1},x_{m-2},\ldots,x_1) = e_{s}$,
$\{s\}$ appears with coefficient $1$.

Since $e_t \theta = e_{t \theta}$,
there exists $\pi$ in the
column stabilizer of $t\theta$ such that $\{s\}=\{t\,\theta\,\pi\}$.
Let $u\in{\widehat{m-1}}$. Then by construction
$x_u\in V_{u+1}(s)\backslash H_u(s)$; since
$\{ s \} = \{ t \theta \pi \}$, it follows that
$x_u\not\in H_u(t\,\theta\,\pi)$.
As $\pi^{-1}$ is a column permutation of $t\theta$, we have
$x_u\in V_{u+1}(t\theta)\cup\ldots\cup V_m(t\theta)$.
Thus
\begin{equation}\label{E:where_is_xu_2}
\forall u \in {\widehat{m-1}},\ \ \ \
x_u\theta^{-1}\in V_{u+1}(t)\cup\ldots\cup V_m(t).
\end{equation}
In particular, $\theta$ does not fix any of the $m-1$
distinct symbols $x_1,\ldots,x_{m-1}\in{\widehat{n-1}}$.

In this paragraph, we will show that $\theta$ does not fix
$n$. Assume that $\theta$ does fix $n$. If the symbols in the
list $y_1, \ldots, y_i$ were distinct, $\theta$ would
be the cycle $(y_i, y_{i-1}, \ldots, y_1, n)$; since
$\theta$ fixes $n$, it follows that there is some repetition
in the  list $y_1, \ldots, y_i$. Since $\theta =
(y_i,n)(y_{i-1},n)\ldots(y_1,n)$ and
$\theta$ fixes $n$,
the only symbols
potentially moved by $\theta$ are on the list
$y_1, \ldots, y_i$. Since this list contains a repeat,
$\theta$ moves at most $i-1$ symbols. The previous paragraph
shows that $\theta$ moves at least $m-1$ symbols. Therefore
$m \leq i$. But by hypothesis $i \leq m-1$. This
contradiction shows that $\theta$ moves $n$.

We now know that $\theta$ moves all the $m$ symbols
in $\{ x_1, \ldots, x_{m-1}, n \}$.
Since $\theta = (y_i,n)(y_{i-1},n)\ldots(y_1,n)$,
$\theta$ can only move symbols on the list
$y_1, y_2, \ldots, y_i, n$. By hypothesis,
$i \leq m-1$. It follows that $i=m-1$, which
is part (a) of our goal. It also follows that
the sets $\{ x_1, \ldots, x_{m-1} \}$ and
$\{ y_1, \ldots, y_{m-1} \}$ coincide and that
the elements on the list $y_1, y_2, \ldots, y_{m-1}$ are
distinct. Hence $\theta$ is equal to the
$m$-cycle $(y_{m-1}, y_{m-2}, \ldots, y_1, n)$.
In particular, $y_{m-1} \theta^{-1} = n$.
  From \eqref{E:where_is_xu_2} applied with $u=m-1$,
$x_{m-1} \theta^{-1} = n$. (This is because $n$ is
the only symbol moved by $\theta$ that is in $V_m(t)$.)
Hence
$y_{m-1} = x_{m-1}$. From this fact and
\eqref{E:where_is_xu_2} applied with $u = m-2$,
it follows that $x_{m-2} \theta^{-1} = x_{m-1}$.
Hence $y_{m-2} = x_{m-2}$. Continuing in this way,
by reverse induction on $u$, it follows that
for all $u \in \widehat{m-1}$, $y_u=x_u$. This gives
goal (b) above, and completes the proof.
\end{proof}

The corresponding result for the action of $L_{n+1}$ on $\ind$ is:

\begin{Lemma}\label{L:L_n+1^i}
Let $t$ be a $\lambda$-tableau and let $T$ be its extension
to $[\lambda{\uparrow^{m+1}}]$. For each $u\in{\widehat m}$,
choose $x_u\in V_u(t)\backslash H_{u-1}(t)$.
Set $S=T\,(n+1,x_m,x_{m-1},\ldots,x_1)$.
Let $i$ be a positive integer with $i \leq m$.
Then the multiplicity of $\{S\}$ in the expansion of
$e_T^\lambda (L_{n+1})^i$ into tabloids is
\[
\begin{aligned}{}
&0,\quad\mbox{when $0\leq i\leq m-1$;}\\
&1,\quad\mbox{when $i=m$.}
\end{aligned}
\]
\end{Lemma}

\begin{proof}
Clearly we have
$(L_{n+1})^i=\sum(w_i,n+1)(w_{i-1},n+1)\ldots(w_1,n+1)$,
where $(w_1,\ldots,w_i)$ ranges over all functions
${\widehat i}\rightarrow{\widehat n}$. Let
$(y_1, \ldots, y_i)$ be a function ${\widehat i}\rightarrow{\widehat n}$,
let
$\theta=(y_i,n+1)(y_{i-1},n+1)\ldots(y_1,n+1)$, and
assume that $\{S\}$ appears with nonzero multiplicity in the expansion
of $e_T^\lambda\theta$ as a linear combination of tabloids.
Then there exists $\pi$ in the column stabilizer of $t\theta$
such that $\{S\}=\{T\,\theta\,\pi\}$.

As $\pi$ fixes the single entry in the last row of $T\theta$, and $x_m$
occupies this node in $S$, it follows that $(n+1)\theta=x_m$. Let
$u\in{\widehat{m-1}}$ and let $s$ denote the restriction of $S$ to
$\lambda$. Then $x_u\in V_{u+1}(s)\backslash H_u(s)$, whence
$x_u\not\in H_u(t\,\theta\,\pi)$. As $\pi^{-1}$ is a column permutation
of $t\theta$, we have $x_u\in V_{u+1}(t\theta)\cup\ldots\cup V_m(t\theta)$.
Thus
\begin{equation}\label{E:where_is_xu_1}
x_u\theta^{-1}\in V_{u+1}(t)\cup\ldots\cup V_m(t).
\end{equation}
In particular, $\theta$ does not fix $x_u$.

  From its definition, $\theta$ moves at most $i+1$ elements of
${\widehat{n+1}}$. But $\theta$ does not fix any of the $m+1$
distinct symbols $n+1,x_m,\ldots,x_1$, and $i\leq m$.
So we must have $i=m$. Together with \eqref{E:where_is_xu_1}, 
this implies
that $x_u\theta^{-1}\in\{x_{u+1},\ldots,x_m\}$.
Reverse induction on $u$ shows that $x_u\theta^{-1}=x_{u+1}$.
Thus $\theta$ coincides with the $(m+1)$-cycle
$(n+1,x_m,x_{m-1},\ldots,x_2,x_1)$. We conclude that
$x_u=y_u$, for $u\in{\widehat m}$. This shows that $\theta$ occurs with
multiplicity $1$ in the expansion of $(L_{n+1})^m$ as a linear
combination of group elements, whence $\{S\}$ appears with multiplicity
$1$ in the expansion of $e_T^{\lambda} (L_{n+1})^m$ as a linear combination
of tabloids in $M^{\lambda{\uparrow^{m+1}}}$. 
\end{proof}

We can now prove the main result of this section.

\begin{Theorem}\label{T:min-poly}
The minimal polynomial of $E_{n-1}$ acting on $\res$ is
\[
\prod_{u=1}^m(x-E(\lambda{\downarrow_u})),
\]
while the minimal polynomial of $E_{n+1}$ acting on $\ind$ is
\[
\prod_{u=1}^{m+1}(x-E(\lambda{\uparrow^u})).
\]
\end{Theorem}

\begin{proof}
First, we will prove the result on $\res$.
Let $t$ be as in Lemma \ref{L:L_n^i}. Then
Lemma \ref{L:L_n^i} implies that the set of vectors
$\{e_t (L_n)^i\mid 0\leq i\leq m-1\}$ is linearly independent.
It follows from Lemma \ref{L:poly} that the set
$\{e_t (E_{n-1})^i\mid 0\leq i\leq m-1\}$ is linearly independent.
So the minimal polynomial of $E_{n-1}$ has degree at least $m$.
But Lemma \ref{L:En-action} and the proof of
Corollary \ref{C:min-poly-upper-bound} show that
$\res  \prod_{u=1}^m(E_{n-1}-E(\lambda{\downarrow_u}))=0$.

The result on $\ind$ follows from an identical argument using
Lemma \ref{L:L_n+1^i} in place of Lemma \ref{L:L_n^i}.
\end{proof}

\section{The indecomposable components of the restriction
          and induction of a Specht module}

The purpose of this section is to compute the indecomposable
components of $\res$ and $\ind$, when the characteristic of
$F$ is not $2$. It is convenient to consider an $F\Sigma_n$-module
$M$ that shares the following properties in common with $\res$ and $\ind$:
\begin{enumerate}
\item
$M$ has a Specht series
\[
0=M_0\subset M_1\subset\ldots\subset M_m=M,
\]
such that $M_u/M_{u-1}\cong S^{\lambda_u}$, where $\lambda_u$
is a partition of $n$, for each $u\in{\widehat m}$.
\item
The labelling partitions satisfy
$\lambda_1\triangleleft\ldots\triangleleft\lambda_m$.
\item
There exists $z\in Z(F\Sigma_n)$ such that the minimal
polynomial of $z$ acting on $M$ has degree $m$.
\end{enumerate}
Looking at the proof of Corollary \ref{C:min-poly-upper-bound},
we see that $z$ has minimal polynomial $\prod_{u=1}^m(x-z_u)$,
where $z$ acts as the scalar $z_u$ on the Specht factor $M_u/M_{u-1}$.

\begin{Lemma}\label{L:tau}
There exists $\tau\in M$ such that 
for all $u\in{\widehat m}$,
$\tau\prod_{i=u+1}^m(z-z_i)$
lies in $M_u\backslash M_{u-1}$.
\end{Lemma}

\begin{proof}
The hypothesis on the degree of the minimal polynomial of $z$
implies that there exists $\tau\in M$ such that $\tau z^{m-1}$
does not lie in the span of the vectors
$\{\tau,\tau z,\ldots,\tau z^{m-2}\}$.
Set $\tau_u=\tau\prod_{i=u+1}^m(z-z_i)$. Repeated application
of Lemma \ref{L:central} shows that $\tau_u\in M_u$.
We claim that $\tau_u\not\in M_{u-1}$. Suppose otherwise.
Then $\tau_u\prod_{i=1}^{u-1}(z-z_u)\subseteq
        M_{u-1}\prod_{i=1}^{u-1}(z-z_u)=0$,
again using Lemma \ref{L:central}.
Thus $\tau\prod_{i=1,i\ne u}^m(z-z_i)=0$.
This contradicts our choice of $\tau$.
\end{proof}

We now consider the endomorphism ring of a module
that has properties (1) and (2) in common with $M$.

\begin{Lemma}\label{L:Specht-series}
Suppose that $\car F \ne 2$. Let $\theta$ be a
$F\Sigma_n$-endomorphism of $M$. Then
\begin{enumerate}
\item
for all $u\in{\widehat m}$, $M_u\theta\subseteq M_u$;

\item
for all $u\in{\widehat m}$,
there is a well-defined $\Sigma_n$-endomorphism
$\theta_u:M_u/M_{u-1}\rightarrow M_u/M_{u-1}$
given by $(v+M_{u-1})\theta_u=v\theta+M_{u-1}$;
\item
the map
$\Phi:\End_{F\Sigma_n}(M)\rightarrow\bigoplus_u
       \End_{F\Sigma_n}(M_u/M_{u-1})$
given by $(\theta) \Phi =(\theta_1,\ldots,\theta_m)$
is an algebra homomorphism;
\item
the kernel of $\Phi$ is the Jacobson
radical of $\End_{F\Sigma_n}(M)$.
\end{enumerate}
\end{Lemma}

\begin{proof}
First, we prove (i).
By induction, we may assume that $M_{u-1}\theta\subseteq M_{u-1}$.
Suppose that $M_u\theta\not\subseteq M_u$.
Choose $v$ so that $m\geq v>u$ and $v$ is
maximal so that $M_u\theta\not\subseteq M_{v-1}$.
Then $M_{u} \theta \subseteq M_{v}$, and 
applying $\theta$ to elements of $M_{u}$ induces
a well-defined nonzero $\Sigma_n$-homomorphism
\[
M_u/M_{u-1}\rightarrow 
M_v/M_{u-1}\twoheadrightarrow M_v/M_{v-1}.
\]
But $\lambda_u\triangleleft\lambda_v$.
This, together with the fact that $\car F\ne 2$,
contradicts 13.17 of \cite{James}, proving (i).
Part (ii)  follows easily from part (i).

It is immediate from the definition of $\theta_u$ that $\Phi$
is an algebra homomorphism. As $\car F\ne 2$, the only
$\Sigma_n$-endomorphisms of $M_u/M_{u-1}$ are scalar multiples
of the identity, by 13.17 of \cite{James}.
It follows that the codomain of $\Phi$ is commutative and semisimple.
Any element of the kernel must send $M_u$ to $M_{u-1}$ for all $u$;
therefore the kernel is nilpotent. 
\end{proof}

We now compute the indecomposable summands of $M$.

\begin{Proposition}\label{P:indM}
Assume that $\car F\neq 2$. Let $B$ be a block idempotent of
$F\Sigma_n$. Then the $F\Sigma_n$-module $MB$ is $0$ or indecomposable.
\end{Proposition}

\begin{proof}
Assume that $MB\neq 0$. Let $A$ be the algebra $\End_{F\Sigma_n}(MB)$.
Identify the algebra $A$ in the natural way with a direct summand of
the algebra $\End_{F\Sigma_n}(M)$. We will use the notation and
results from Lemma \ref{L:Specht-series} throughout this proof.
Our goal is to show that $A/J(A)$ has dimension $1$ over $F$.

Suppose then that $\theta\in A$.
Let $w$ be 
maximal such that the Specht module $M_w/M_{w-1}$
belongs to $B$. Our task is to show that if $\theta_w=0$,
then $\theta_u=0$ for all $u$ such that $M_u/M_{u-1}$ belongs to $B$.
(The proposition follows easily from this. Let $\phi$ be in $A$.
Then there is a scalar $c$ such that
the map $\phi_w$ is $c$ times the
identity. Let $\theta = \phi - c 1_{A}$.
Then $\theta_w = 0$. Since $\theta_{u}$ is also $0$ for all
$u$ with $M_u/M_{u-1}$ belonging to $B$, it follows
from the last part of Lemma \ref{L:Specht-series}
that $\theta \in J(A)$. Hence $A/J(A)$ has dimension $1$.)

Now assume that $\theta_w=0$, and let $u$ be an integer
such that $M_u/M_{u-1}$ belongs to $B$.
Let $\tau\in M$ be as in  Lemma \ref{L:tau},
set  $\tau_{u}: = \tau \prod_{i=u+1}^{m}(z-z_{i})$,
and set $\tau_{w}: = \tau \prod_{i=w+1}^{m}(z-z_{i})$.
The lemma states that $\tau_u\in M_u\backslash M_{u-1}$
and $\tau_w\in M_w\backslash M_{w-1}$.
Since $u \leq w$,
we have
\[
\begin{aligned}{}
\tau_u\theta
          &=\left(\tau_w\prod_{i=u+1}^w(z-z_i)\right)\theta\\
          &=\tau_w\theta\prod_{i=u+1}^w(z-z_i),
                  \quad\mbox{as $\theta$ is in 
          $\End_{F\Sigma_n}(M)$,}\\
          &\in M_{w-1}\prod_{i=u+1}^w(z-z_i),
\quad\mbox{as $\theta_w=0$ implies that $\tau_w\theta\in M_{w-1}$,}\\
          &=\left(M_{w-1}\prod_{i=u+1}^{w-1}(z-z_i)\right)(z-z_w)\\
          &\subseteq M_u(z-z_w),
                  \quad\mbox{using  Lemma \ref{L:central}
                                  repeatedly}.
\end{aligned}
\]
Now $M_u/M_{u-1}$ and $M_w/M_{w-1}$ both belong to $B$. So $z_u=z_w$, 
since both scalars are equal to 
the image of $z$ under the central character of $B$. Lemma
\ref{L:central} and the last inclusion displayed
above then show that $\tau_u\theta\in M_{u-1}$.
But $\tau_u\not\in M_{u-1}$, as proved in Lemma \ref{L:tau},
and $\End_{F\Sigma_n}(M_u/M_{u-1})$ is one-dimensional, by 13.17
of \cite{James}. We conclude that $\theta_u=0$, as required. 
\end{proof}

We have now done all the work to prove the main result of this paper.

\begin{Theorem}
\label{Main}
Assume that $\car F\neq 2$.
Let $b$ be a block idempotent of $F\Sigma_{n-1}$.
Then the $F\Sigma_{n-1}$-module
$(S^\lambda{\downarrow}_{{S_{n-1}}})\,b$
is $0$ or indecomposable. Let $B$ be a block idempotent
of $F\Sigma_{n+1}$. Then the $F\Sigma_{n+1}$-module
$(S^\lambda{\uparrow}^{{S_{n+1}}})\,B$
is $0$ or indecomposable.
\end{Theorem}

\begin{proof}
We know that $\res$ and $\ind$ satisfy properties (1) and (2) of $M$.
That they also satisfy property (3) is a consequence of Theorem
\ref{T:min-poly}. The result now follows from Proposition \ref{P:indM}.
\end{proof} 

We will finish by giving examples to show that the 
assumption $\car F \neq 2$ cannot be dropped in Theorem
\ref{Main}. 

Assume that $\car F = 2$. Consider the Specht module 
$S^{(6,1,1,1)}$. The decomposition matrix for $\Sigma_{9}$ given
in \cite{James} shows that $S^{(8,1)}$ and $S^{(6,3)}$ are simple
and that $S^{(6,1,1,1)}$ has a composition series with factors
$S^{(8,1)}$ and $S^{(6,3)}$. By 23.8 in \cite{James}, 
$S^{(6,1,1,1)}$ is self-dual, so there is another composition series
in which the factors appear in the 
other order. It follows that $S^{(6,1,1,1)} \cong S^{(8,1)} \oplus 
S^{(6,3)}$. 

Now consider the restriction of $S^{(6,1,1,1)}$ to $\Sigma_{8}$. 
All components of the restriction belong to the
principal $2$-block of $\Sigma_{8}$, which is the 
block with empty core.
Since $S^{(6,1,1,1)}$ is decomposable, 
so is its restriction to $\Sigma_{8}$.

For the other counterexample, 
let $M = S^{(6,1,1)}\uparrow ^{\Sigma_{9}}$. The module
$M$ has a Specht series with factors $S^{(7,1,1)}$, 
$S^{(6,2,1)}$, and $S^{(6,1,1,1)}$. These factors belong
to 2-blocks with cores $(1)$, $(1)$, and $(2,1)$ respectively. 
It follows that if $B$ is the block idempotent corresponding to 
$2$-core $(2,1)$, then $MB \cong S^{(6,1,1,1)}$;
thus $MB$ is decomposable.    

\section{Acknowledgement}

Part of this paper was written while the first author was visiting the
National University of Ireland, Maynooth. The visit was funded by a grant
from Enterprise Ireland, under the International Collaboration
Programme 2003. Enterprise Ireland support is funded under the National
Development Plan and co-funded by European Union Structural Funds.
We gratefully acknowledge this assistance.

Although they now require no computer calculations,
the examples at the end of section 3 were originally found
using computer programs written in GAP and Magma. The
programs were written by Julia Dragan-Chirila, under the
supervision of the first author. Her work
was supported by Northern Illinois University's
Undergraduate Research Apprenticeship Program.

\end{document}